\documentclass{article}
\usepackage[utf8]{inputenc}
\usepackage{relsize}
\usepackage{amsfonts}
\usepackage{amsthm}
\newcommand{\no}{\noindent}
\title{Free dihedral actions on Abelian varieties}
\author{Bruno Aguil\'o Vidal}
\date{}
\begin{document}
\maketitle

\begin{abstract}
\noindent We give a simple construction for hyperelliptic varieties defined as the quotient of a complex torus by the action of a dihedral group that contains no translations and fixes no points. This generalizes a construction given by Catanese and Demleitner for $D_4$ in dimension three.
    
\end{abstract}
\vspace{0.5 cm}

\begin{center}
    \Large{\textbf{Introduction}}
\end{center}

A Generalized Hyperelliptic Manifold $X$ is defined to be a quotient $X = T/G$ of a complex torus $T$ by the free action of a finite group $G$ which contains no translations. We say that X is a Generalized Hyperelliptic Variety if moreover the torus T is projective, i.e., it is an Abelian variety $A$.

Uchida and Yoshihara showed that the only non Abelian group that gives such an action in dimension three is the dihedral group $D_4$ of order $8$ \cite{UchYo}. Later, Catanese and Demleitner gave a simple and explicit construction for that action \cite{CatDem} and completed the characterization of three-dimensional hyperelliptic manifolds \cite{CatDem2}.

The purpose of this note is to generalize Catanese and Demleitner's construction to bigger dihedral groups acting in higher dimension. Specifically, for every $n\in\mathbb{N}$ we give a Generalized Hyperelliptic Variety of dimension $2n+1$ defined by the action of the dihedral group $D_{4n}$ of order $8n$ acting on a family of Abelian varieties, from which the construction by Catanese and Demleitner remains as the particular case for $n=1$. 

We end with a simple corollary that explains how this allows us to create Generalized Hyperellyptic Varieties using any dihedral group.
\vspace{0.5 cm}
\begin{center}
    \Large{\textbf{The construction}}
\end{center}
\noindent Let $E,\,E'$ be any two elliptic curves,
\begin{center}
    $E=\mathbb{C}/\left( \mathbb{Z} + \mathbb{Z}\tau\right)$, $E'=\mathbb{C}/\left( \mathbb{Z} + \mathbb{Z}\tau'\right).$
\end{center}
\vskip0.3cm
\noindent Now, for $n\in\mathbb{N}$ set $A':= E^{2n}\times E'$, $A:= A'/\langle w \rangle$, where $w:= \left(1/2,\, 1/2,\, ...,\, 1/2,\, 0\right)$.

\vskip0.3cm
\noindent \textbf{Theorem.} \textit{The Abelian Variety $A$ admits a free action with no translations of the dihedral group $D_{4n}$ of order $8n$}.

\vskip0.5cm
\begin{proof} First, let us recall that for $k\in\mathbb{N}$, the dihedral group of order $2k$ is defined as 
\begin{center}
    $D_{k} := \left\langle r,s\,|\, r^k = 1,\,s^2 = 1,\, (rs)^2 = 1\right\rangle.$
\end{center}
Now, set, for $z:= \left(z_1,\, z_2,\, ...,\, z_{2n},\, z_{2n+1}\right)\in A'$:

\begin{eqnarray*}
    r(z) & := & \left(-z_{2n},\, z_1,\, z_2,\, ...,\, z_{2n-1},\, z_{2n+1} + \mathsmaller{\frac{1}{4n}}\right)\\ 
    & = & R(z) + \left(0,\,...,\,0,\, \mathsmaller{\frac{1}{4n}}\right), \\
    s(z) & := & \left(-z_{2n} + b_1,\, -z_{2n-1} + b_2,\,...,\, -z_1 + b_{2n},\, -z_{2n+1}\right)\\ 
    & = & S(z) + (b_1,\,b_2,\,...,\,b_{2n},\,0),
\end{eqnarray*}\\
    where, for $i=1,\,...,\,n$, $b_{2i-1} := 1/2 + \tau/2$ and $b_{2i} := \tau/2$.
\vskip0.3cm
\no\textbf{Step 1.} It is easy to verify that $r$ and $R$ have order exactly $4n$ on $A'$, and that $R(w) = w$, so that $r$ descends to an automorphism of $A$ of order exactly $4n$. Moreover, any power $r^j$, $0< j<4n$, acts freely on $A$ since the $(2n+1)$-th coordinate of $r^j (z)$ equals $z_{2n+1} + \frac{j}{4n}$, and clearly none of this powers is a translation.\\

\vskip0.1cm
\no\textbf{Step 2.} $s^2 (z) = z+w$, since for $i = 1,\,...,\, 2n$, $b_i - b_{2n + 1-i}= 1/2$; moreover, $S(w)=w$, hence $s$ descends to an automorphism of $A$ of order exactly $2$.\\
\vskip0.1cm
\no\textbf{Step 3.} We have
\begin{eqnarray*}
    rs(z) & = & r\left(-z_{2n} + b_1,\, -z_{2n-1} + b_2,\,...,\, -z_1 + b_{2n},\,-z_{2n+1}\right)\\
    & = & \left(z_1 - b_{2n},\, -z_{2n} + b_1,\,...,\, -z_2 + b_{2n-1},\,-z_{2n+1} + \mathsmaller{\frac{1}{4n}}\right).
\end{eqnarray*}\\
\no hence\\
\begin{eqnarray*}
    (rs)^2 (z) & = & \left(z_1 -2b_{2n},\,z_2 + b_1 - b_{2n-1},\,...,\, z_i + b_{i-1} - b_{2n - (i-1)},\,...,\,z_{2n} + b_{2n-1} - b_1,\, z_{2n+1}\right)\\
    & = & z,
\end{eqnarray*}\\
and we have an action of $D_{4n}$ on $A$, since the orders of $r$, $s$ and $rs$ are precisely $4n$, $2$ and $2$.\\
\vskip0.1cm
\no\textbf{Step 4.} We claim that also the symmetries in $D_{4n}$ act freely on $A$ and are not translations. Since there are exactly two conjugacy classes of symmetries, those of $s$ and $rs$, it suffices to observe that these two transformations are not translations. In the next step we show that they both act freely.\\
\vskip0.1cm
\no\textbf{Step 5.} It is rather immediate that $rs$ acts freely, since $rs(z) = z$ in $A$ is equivalent to 
\begin{center}
    $\left( -b_{2n},\, -z_{2n} - z_2 + b_1,\,...,\, -z_2 - z_{2n} + b_{2n-1},\, -2z_{2n+1} + \mathsmaller{\frac{1}{4n}}\right)$
\end{center}
being a multiple of $w$ in $A'$, but this is absurd since $2w = 0$ and $-b_{2n} = \tau/2 \neq 0,\,1/2$.\\

On the other hand, $s$ acts freely in $A$ because $s(z) = z$ is equivalent to
\begin{center}
    $\left( -z_{2n} - z_1 + b_1,\,-z_{2n-1}-z_2 + b_2,\,...,\, -z_1 - z_{2n} + b_{2n},\,-2z_{2n+1}\right)$
\end{center}
being a multiple of $w$ in $A'$, but the first and $2n$-th coordinate of multiples of $w$ are equal, while here the difference between them is $1/2\neq 0$.\\
\end{proof}
\vskip0.5cm
\begin{center}
    \Large{\textbf{Using any dihedral group}}
\end{center}

Notice that, although the previous  construction is somewhat restrictive because it works with very specific dihedral groups, since it is true that $D_n\subseteq D_{nk}$ for all $n,\,k\in\mathbb{N}$, we have the following corollary:\\

\vskip0.3cm
\no\textbf{Corollary.} \textit{For all $n\in\mathbb{N}$, there exists a free action of the dihedral group $D_n$ of order $2n$ on some Abelian variety of dimension $\frac{\textbf{mcm}(4,n)}{2} + 1$ that contains no translations}.

\bigskip

\noindent\textit{Acknowledgements}: I would like to thank professors Robert Auffarth and Giancarlo Lucchini Arteche for introducing me to this topic and for dedicatedly guiding me into getting this result as an undergraduate student.

\bigskip

\bibliography{Referencias}
\bibliographystyle{plain}
\vspace{0.5cm}

B. Aguil\'o Vidal, Departamento de Matem\'aticas, Facultad de Ciencias, Universidad de Chile, Las Palmeras 3425, \~{N}u\~{n}oa, Santiago, Chile.\\

\textit{E-mail address}: \textbf{bruno.aguilo@ug.uchile.cl}
\end{document}